\documentclass[12pt,a4paper,twoside]{article}
\usepackage{mathtext}
\usepackage[cp1251]{inputenc}
\usepackage[T2A]{fontenc}
\usepackage[russian]{babel}
\usepackage[dvips]{graphicx}
\usepackage{amsmath}
\usepackage{amssymb}
\usepackage{amsxtra}
\usepackage{latexsym}
\usepackage{ifthen}
\usepackage[centerlast,footnotesize,sl]{caption2}
\begin{document}

\newcounter{lemma}[section]
\newcommand{\lemma}{\par \refstepcounter{lemma}%
{\bf Лемма \arabic{section}.\arabic{lemma}.}}
\renewcommand{\thelemma}{\thesection.\arabic{lemma}}

\newcounter{corollary}[section]
\newcommand{\corollary}{\par \refstepcounter{corollary}%
{\bf Следствие \arabic{section}.\arabic{corollary}.}}
\renewcommand{\thecorollary}{\thesection.\arabic{corollary}}

\newcounter{remark}[section]
\newcommand{\remark}{\par \refstepcounter{remark}%
{\bf Замечание \arabic{section}.\arabic{remark}.}}
\renewcommand{\theremark}{\thesection.\arabic{remark}}

\newcounter{theorem}[section]
\newcommand{\theorem}{\par \refstepcounter{theorem}%
{\bf Теорема \arabic{section}.\arabic{theorem}.}}
\renewcommand{\thetheorem}{\thesection.\arabic{theorem}}

\newcounter{proposition}[section]
\newcommand{\proposition}{\par \refstepcounter{proposition}%
{\bf Предложение \arabic{section}.\arabic{proposition}.}}
\renewcommand{\theproposition}{\thesection.\arabic{proposition}}

\newcommand{\proof}{{\it Доказательство.\,\,}}

\renewcommand{\theequation}{\arabic{section}.\arabic{equation}}

\def\Xint#1{\mathchoice
   {\XXint\displaystyle\textstyle{#1}}
   {\XXint\textstyle\scriptstyle{#1}}
   {\XXint\scriptstyle\scriptscriptstyle{#1}}
   {\XXint\scriptscriptstyle\scriptscriptstyle{#1}}
   \!\int}
\def\XXint#1#2#3{{\setbox0=\hbox{$#1{#2#3}{\int}$}
     \vcenter{\hbox{$#2#3$}}\kern-.5\wd0}}
\def\dashint{\Xint-}

\noindent
\renewcommand{\baselinestretch}{1.25}
\normalsize \large

\noindent{\bf УДК 517.5}\hfill{\bf М А Т Е М А Т И К А}

\medskip

\centerline{\bf А.С. Ефимушкин, В.И. Рязанов}

\medskip

\centerline{\bf О задаче Римана-Гильберта для аналитических функций
} \centerline{\bf в круговых областях}

\medskip

\centerline{\it (Представлено член-корреспондентом НАН Украины В.Я.
Гутлянским)}

\medskip

{\bf 1. Введение.} Краевые задачи для аналитических функций $f$
восходят к знаменитой диссертации Римана (1851), а также известным
работам Гильберта (1904, 1912, 1924), и Пуанкаре (1910), смотри
также историю вопроса в монографии \cite{Vek} и нашей предыдущей
статье  \cite{UMV}.

В 1904 году Гильберт сформулировал следующую задачу, которую принято
называть проблемой Гильберта или {\it проблемой
Римана-Гиль\-бер\-та}. Она состояла в доказательстве существования и
нахождении аналитической функции $f$ в области $D\subset\mathbb{C}$,
ограниченной спрямляемой жордановой кривой $K$ с условием
$$\lim_{z\rightarrow\zeta} \mathrm{Re}\ \overline{\lambda(\zeta)}\cdot f(z)=\varphi(\zeta) \qquad
\forall\ \zeta\in K, \eqno(1) $$ где им предполагалось, что функции
$\lambda$ и $\varphi$ непрерывно дифференцируемы относительно
натурального параметра длины на кривой $K$, и что $|\lambda|\neq0$
на $K$. Поэтому можно считать, что $|\lambda(\zeta)|\equiv 1$.

\medskip

{\bf 2. Определения и предварительные замечания.} За определением и
основными свойствами логарифмической \"ем\-ко\-сти мы отсылаем
читателя к хорошо известным работам \cite{Kar}--\cite{F}, а также к
статье \cite{UMV}.

\medskip

Напомним, см. \cite{Go}, что область $D$ в
$\overline{\mathbb{C}}=\mathbb{C}\cup\{\infty\}$ именуется {\it
круговой}, если ее граница состоит из конечного числа взаимно
непересекающихся окружностей и точек. Если граница этой области
состоит только из окружностей, то называем эту область {\it
невырожденной}.

Пусть $\mathbb{D}_*$ -- ограниченная невырожденная круговая область
в $\mathbb{C}$. Мы называем
$\lambda:\partial\mathbb{D}_*\to\mathbb{C}$ {\it функцией
ограниченной вариации} на открытой дуге
$\gamma\subset\partial\mathbb{D}_*$, пишем $\lambda~\in
\mathcal{BV}(\gamma)$, если

$$V_{\lambda}(\gamma)\ :=\ \sup\ \sum_{j=0}^{k-1}\
|\lambda(\zeta_{j+1})-\lambda(\zeta_j)|<\infty\ ,   \eqno(2)$$ где
супремум бер\"ется по всем конечным наборам точек $\zeta_j
\in\gamma$, $j=0,1,\dots,k$, таким, что $\zeta_j$ лежит между
$\zeta_{j+1}$ и $\zeta_{j-1}$ для каждого $j=1,\dots,k-1$.

Мы также называем $\lambda$ {\it функцией счетно-ограниченной
вариации}, пишем $\lambda\in \mathcal{CBV}(\partial\mathbb{D}_*)$,
если найдется счетное число попарно-непересекающихся открытых дуг
$\gamma_n\subset\partial\mathbb{D}_*$, на которых она является
функцией ограниченной вариации, а множество
$\partial\mathbb{D}_*\setminus\bigcup\gamma_n$ счетно.

\medskip

Рассуждая  как в случае класса $\mathcal{BV}(\partial\mathbb{D})$,
где $\mathbb{D}=\{z\in\mathbb{C}:|z|~<1\}$, полностью аналогично
доказательству предложения 5.1 в нашей работе \cite{UMV}, получаем
его аналог.

{\bf Предложение 1. } {\it  Для любой функции
$\lambda:\partial\mathbb{D}_*\to\partial\mathbb{D}$ класса
$\mathcal{CBV}(\partial\mathbb{D}_*)$ найдется функция
$\alpha_{\lambda}: \partial\mathbb{D}_*\to(-\pi,\pi]$ класса
$\mathcal{CBV}(\partial\mathbb{D}_*)$, такая, что
$\lambda(\zeta)=\exp\{{i\alpha_{\lambda}}(\zeta)\}$,
$\zeta\in\partial\mathbb{D}_*$.}

\medskip

В дальнейшем мы называем функцию $\alpha_{\lambda}$ {\it функцией
аргумента} $\lambda$.

\medskip

По теореме 1 в \cite{T}, рассуждая аналогично секции 5 в \cite{UMV},
получаем также аналог теоремы 5.1.

\medskip

{\bf Лемма 1. } {\it Пусть $\alpha:\partial\mathbb D\to\mathbb R$ -
ограниченная функция класса $\mathcal{CBV}(\partial\mathbb{D})$ и
пусть $f:\mathbb D\to\mathbb C$ - аналитическая функция, такая, что
$$
\lim\limits_{z\to\zeta}\ \mathrm {Re}\ f(z)\ =\ \alpha(\zeta)
\quad\quad\quad\mbox{для}\ \ \mbox{п.в.}\ \ \
\zeta\in\partial\mathbb D \eqno(3)$$ относительно логарифмической
\"емкости вдоль любых некасательных путей. Тогда
$$
\lim\limits_{z\to\zeta}\ \mathrm {Im}\ f(z)\ =\ \beta(\zeta)
\quad\quad\quad\mbox{для}\ \ \mbox{п.в.}\ \ \
\zeta\in\partial\mathbb D \eqno(4)$$ относительно логарифмической
\"емкости вдоль любых некасательных путей, где
$\beta:\partial\mathbb D\to\mathbb R$ - некоторая функция измеримая
относительно логарифмической \"емкости. }

\medskip

{\bf 3. Основные результаты.} Следующие теоремы дают решение задачи
Римана-Гильберта в круге и круговых областях.

\medskip

{\bf Теорема 1. } {\it Пусть $\lambda:\partial\mathbb
D\to\partial\mathbb D$ - функция  класса
$\mathcal{CBV}(\partial\mathbb{D})$ и $\varphi:\partial\mathbb
D\to\mathbb R$ - функция измеримая относительно логарифмической
\"емкости. Тогда существует аналитическая функция $f:\mathbb
D\to\mathbb C$, такая, что вдоль любых некасательных путей
$$
\lim\limits_{z\to\zeta}\ \mathrm {Re}\
\{\overline{\lambda(\zeta)}\cdot f(z)\}\ =\ \varphi(\zeta)
\quad\quad\quad\mbox{для}\ \ \mbox{п.в.}\ \ \
\zeta\in\partial\mathbb D \eqno(5)$$ относительно логарифмической
\"емкости.}

\medskip

Дуйствительно,  по предложению 1 функция аргумента
$\alpha_{\lambda}$ принадлежит классам $L^{\infty}(\partial\mathbb
D)$ и $\mathcal{CBV}(\partial\mathbb{D})$. Поэтому
$$ g(z)\ =\ \frac{1}{2\pi i}\ \int\limits_{\partial\mathbb
D}\alpha(\zeta)\ \frac{z+\zeta}{z-\zeta}\
 \frac{d\zeta}{\zeta}\ , \ \ \ \ \ z\in\mathbb D\ ,  \eqno(6)
$$
является аналитической функцией в $\mathbb D$ с $u(z):={\mathrm Re}\
g(z)\to\alpha(\zeta)$ при $z\to\zeta$ вдоль любых путей в $\mathbb
D$ для всех точек $\zeta\in\partial\mathbb D$ за исключением их
счетного числа, см., например, теорему IX.1.1 в \cite{Go}, с.369, и
теорему I.D.2.2 в \cite{Ku}, с. 18. Отметим, что ${\cal
A}(z)=\exp\{ig(z)\}$ также является аналитической функцией.

По лемме 1 существует функция $\beta:\partial\mathbb D\to\mathbb R$
конечная п.в. и измеримая относительно логарифмической \"емкости,
такая, что $v(z):={\mathrm Im}\ g(z)\to\beta(\zeta)$  вдоль любых
некасательных путей при $z\to\zeta$ для п.в.
$\zeta\in\partial\mathbb D$ также относительно логарифмической
\"емкости. Таким образом, по следствию 4.1 в \cite{UMV} существует
аналитическая функция ${\cal B}:\mathbb D\to\mathbb C$, такая, что
$U(z):=\mathrm {Re}\ {\cal B}(z)\to\varphi(\zeta)\cdot
\exp\{{\beta(\zeta)}\}$ при $z\to\zeta$ вдоль любых некасательных
путей для п.в. $\zeta\in\partial\mathbb D$ относительно
логарифмической \"емкости. Наконец, элементарные вычисления
показывают, что искомая функция $f={\cal A}\cdot{\cal B}$.

\medskip

{\bf Теорема 2. } {\it Пусть $\mathbb{D}_*$ -- ограниченная
невырожденная круговая область,
$\lambda:\partial\mathbb{D}_*\to\partial\mathbb{D}$~-- функция
класса $\mathcal{CBV}(\partial\mathbb{D}_*)$ и
$\varphi:\partial\mathbb{D}_*\to\mathbb{R}$~-- функция измеримая
относительно логарифмической ёмкости. Тогда существует многозначная
аналитическая функция $f:\mathbb{D}_*\to\mathbb{C}$, такая, что
вдоль любых некасательных путей
$$
\lim_{z\to\zeta}\mathrm{Re}\{\overline{\lambda(\zeta)}\cdot
f(z)\}=\varphi(\zeta)\quad\quad\quad\mbox{для}\ \ \mbox{п.в.}\ \ \
\zeta\in\partial\mathbb {D}_* \eqno(7)
$$
относительно логарифмической ёмкости}.

\medskip

Действительно, по теореме Пуанкаре, см., например, теорему 6.1 в
\cite{Go}, существует локально конформное отображение $g$ круга
$\mathbb{D}$ на $\mathbb{D}_*$. Пусть
$h:\mathbb{D}_*\to\mathbb{D}$~-- соответствующая многозначная
аналитическая функция, которая является обратной к $g$.

Область $\mathbb{D}_*$ без конечного числа разрезов является
односвязной областью и, следовательно, $h$ имеет там однозначные
ветви, продолжимые на границу по теоремам Каратеодори. Согласно
разделу VI.2 в \cite{Go}, $\partial\mathbb{D}$ за исключением
счетного множества своих точек, состоит из счетной совокупности дуг,
каждая из которых взаимно однозначно (гомеоморфно) отображается $g$
на некоторую окружность из $\partial\mathbb{D}_*$ с одной
выброшенной точкой.

Отметим, что по принципу отражения $g$ конформно продолжима в
окрестность каждой из таких дуг и, следовательно, некасательные пути
к точкам этих дуг из $\partial\mathbb{D}$ переходят в некасательные
пути к соответствующим точкам на окружностях из
$\partial\mathbb{D}_*$ и наоборот.

Пусть $\Lambda=\lambda\circ g$ и $\Phi=\varphi\circ g$, где $g$
считаем уже продолженным на указанные дуги единичной окружности
$\partial\mathbb{D}$. Заметим, что из условий теоремы по построению
следует, что $\Lambda\in\mathcal{CBV}(\partial\mathbb{D})$, а
функция $\Phi=\varphi\circ g$ измерима относительно логарифмической
емкости.

По теореме 1 существует аналитическая функция
$F:\mathbb{D}\to\mathbb{C}$, такая, что вдоль любых некасательных
путей
$$\lim_{w\to\eta}\mathrm{Re}\{\overline{\Lambda(\eta)}\cdot F(w)\}=\Phi(\eta)\quad\quad\quad\mbox{для}\ \ \mbox{п.в.}\ \ \
\eta\in\partial\mathbb D \eqno(8)$$ относительно логарифмической
ёмкости.

В силу приведенных выше рассуждений, $f:=F\circ h$ является искомой
многозначной аналитической функцией, удовлетворяющей граничному
условию (7).

{\bf Теорема 3. } {\it В условиях теорем 1 и 2 пространства решений
задачи Римана-Гильберта имеют бесконечную размерность. }

Последний результат является прямым следствием конструкций решений в
теоремах 1 и 2, а также теоремы 2.1 из работы \cite{R2}.

\newpage
{\bf Про задачу Рімана-Гільберта для аналітичних функцій у кругових
областях. --} {\bf Єфімушкін А.С., Рязанов В.І.}

\medskip

Доведено існування однозначних аналітичних розв'язків в одиничному
колі та багатозначних аналітичних розв'язків в областях, обмежених
скінченним числом кіл, задачі Рімана-Гільберта із коефіцієнтами
злічено-обмеженної варіації та граничними даними, що є вимірюваними
відносно логарифмічної ємності. Показано, що простори розв'язків
мають нескінчену розмірність.

{\bf Ключовi слова:} задача Рімана-Гільберта, аналітичні функції,
кругові області.

\medskip

{\bf О задаче Римана-Гильберта для аналитических функций в круговых
областях. --} {\bf Ефимушкин А.С., Рязанов В.И.}

\medskip

Доказано существование однозначных аналитических решений в единичном
круге и многозначных аналитических решений в областях, ограниченных
конечным числом окружностей, задачи Римана-Гильберта с
коеффициентами счетно-ограниченной вариации и граничными данными,
измеримыми относительно логарифмической емкости. Показано, что
пространства решений имеют бесконечную размерность.

{\bf Ключевые слова:} задача Римана-Гильберта, аналитические
функции, круговые области.

\medskip

{\bf On the Riemann-Hilbert problem for analytic functions in
circular domains. --} {\bf Yefimushkin A.S., Ryazanov V.I.}

\medskip

It is proved the existence of single-valued analytic solutions in
the unit disk and multivalent analytic solutions in domains bounded
by a finite collection of circles for the Riemann-Hilbert problem
with coefficients of sigma finite variation and with boundary data
that are measurable with respect to logarithmic capacity. It is
shown that these spaces of solutions have the infinite dimension.


\newpage

Ефимушкин А.С.,

Институт математики НАН Украины,

ул. Терещенковская 3, 01601, Киев,

(050)-049-00-74, a.yefimushkin@gmail.com

\bigskip

Рязанов Владимир Ильич,

Институт прикладной математики

и механики НАН Украины,

ул. Добровольского 1, 84100, Славянск,

(050)-584-90-65, vl.ryazanov1@gmail.com

\end{document}